\documentclass{elsart1}

\usepackage{graphicx}

\usepackage{amssymb}

\usepackage{amsfonts}
\usepackage{amsmath}

\newcommand{\N}{{\mathbb N}}
\newcommand{\Z}{{\mathbb Z}}
\newcommand{\R}{{\mathbb R}}

\begin{document}

\begin{frontmatter}

\title{On computing B\'ezier curves by Pascal matrix methods}

\author[sul]{Licio Hernanes Bezerra\corauthref{cor}}
\corauth[cor]{Corresponding author.}
\ead{licio@mtm.ufsc.br}
\author[norte]{Leonardo Koller Sacht}
\ead{leo-ks@impa.br}

\address[sul]{Universidade Federal de Santa Catarina, Departamento de Matem\'atica,
Florian\'opolis, SC 88040-900, Brazil}
\address[norte]{Instituto de Matem\'atica Pura e Aplicada,
Rio de Janeiro, RJ 22460-320, Brazil}

\begin{abstract}
The main goal of the paper
is to introduce methods which compute B\'ezier curves faster than Casteljau's method does.
These methods are based on the spectral factorization of a $n\times n$ Bernstein matrix,
$B^e_n(s)= P_nG_n(s)P_n^{-1}$,
where $P_n$ is the $n\times n$ lower triangular Pascal matrix.
So we first calculate the exact optimum positive value $t$ in order
to transform $P_n$ in a scaled Toeplitz matrix, which is a problem
that was partially solved by X. Wang and J. Zhou (2006).
Then fast Pascal matrix-vector multiplications and strategies of polynomial evaluation
are put together to compute B\'ezier curves.
Nevertheless, when $n$ increases, more precise Pascal matrix-vector
multiplications allied to affine transformations of the vectors of coordinates of the control points
of the curve
are then necessary to stabilize all the computation.
\end{abstract}

\begin{keyword}
Pascal matrix \sep Bernstein polynomial \sep B\'ezier curve  \sep Toeplitz matrix
\MSC 15A18  \sep 65F15 \sep 68U07
\end{keyword}
\end{frontmatter}

\section{Introduction}

B\'ezier has his name on the curve $B$ of degree $n-1$ defined from  $n$ given points
$Z_0 = (x_0,y_0)$, $Z_1=(x_1,y_1)$, ..., $Z_{n-1}=(x_{n-1},y_{n-1})$
in $\mathbb{R}^2$
as follows:
$$
B(s) = (x(s),y(s)) = \sum_{i=0}^{n-1} Z_ib_{i,n-1}(s), \hspace{1cm} s \in [0,1],
$$
where $b_{i,n-1}(s) = \binom{n-1}{i} s^i (1-s)^{n-1-i}$, for each $i\in \{ 0,...,n-1\}$, which is called
a Bernstein polynomial \cite{bezier}.
That is, for each $s \in [0,1]$ we have
$$
x(s) = \sum_{i=0}^{n-1} b_{i,n-1}(s)x_i, \quad y(s) = \sum_{i=0}^{n-1}b_{i,n-1}(s)y_i.
$$
Notice that $B(s) = B_r(1-s)$, where $B_r(s)$ is the B\'ezier
curve defined from the control points in the reverse order:
$$B_r(s) = \sum_{i=0}^{n-1} Z_{n-1-i}b_{i,n-1}(s).$$
Let $s$ be a real number. The $n\times n$ lower triangular matrix $B^e_n(s)$ such that
$\left( B^e_n\right)_{ij}(s) = b_{j-1,i-1}(s)$, for each $n\ge i\ge j\ge 1$, is called a Bernstein matrix \cite{aceto}.

Paul de Casteljau developed a very stable algorithm to evaluate B\'ezier curves.
In this so called Casteljau's algorithm
the number of arithmetic operations grows quadratically with $n$:
it requires $n(n-1)/2$ additions and $n(n-1)$ multiplications to
calculate a point (which is not an endpoint) on a B\'ezier curve of degree $n-1$
\cite{casteljau}.
A natural question is if there could be a less expensive algorithm to
compute a B\'ezier curve, which is answered e.g. in \cite{phien} for $n < 10$.
Here, we extend this answer to $n \le 64$
by a
different approach, which  arises from the
expression of a $n\times n$ Bernstein matrix in terms of the lower triangular
Pascal matrix $P_n$ (see \cite{aceto}) as follows:
$$
x(s) = e_n^TP_nG(-s)P_nG(-1) x,
\quad y(s) = e_n^TP_nG(-s)P_nG(-1) y.
$$
Here, $x= (x_0 \, x_1 \, \cdots x_{n-1})^T$,
$y= (y_0 \, y_1 \, \cdots y_{n-1})^T$;
$e_n$ denotes the nth canonical vector of $\R^n$; $G(s)$ is the diagonal matrix
such that $e_k^T G(s) e_k = s^{k-1}$ for all $k\in \{ 1,...,n\}$; and $P_n$ is defined
by
$$\left( P_n\right)_{ij}= \left\{ \begin{array}{ccl}
\binom{i-1}{j-1} &,& \mbox{for } i \geqslant j; \\
0&,& \mbox{otherwise}.
\end{array}
\right.$$
From now on, the methods introduced here which are
originated from the above expression will be called
Pascal matrix methods.
One central problem in these methods is how we can
accurately compute a matrix-vector multiplication with the $n\times n$ lower triangular Pascal matrix.
In \S 2, we present an algorithm that utilizes only n(n-1)/2 additions
that computes a matrix-vector multiplication with $P_n$ in a very precise way, which is based on
the property of the lower triangular Pascal matrix be a product
of bidiagonal matrices.
But this time complexity can be reduced to ${\cal O}(n\log\, n)$
by using the fact that the lower triangular Pascal matrix
is similar to lower triangular Toeplitz matrices via diagonal
matrices \cite{wang}. This transformation depends on a
positive real parameter $t$, which is arbitrary. In \cite{wang}, it is suggested that $(n-1)/e$ could be taken
as a good approximate value for the optimum parameter. Here, the optimum value when it exists is calculated.
It is left as an open problem if the set of $n$ for which it does not exist the optimum value
is finite.
At the end of the section, there are results of tests to compare the accuracy of
matrix-vector multiplications carried out by the fast multiplication algorithm,
with both approximate and exact optimum values, and by
the algorithm based on the bidiagonal factorization.

The other central problem in Pascal matrix methods is, once calculated $z=P_nG(-1)x$ (or $z=P_nG(-1)y$),
how to efficiently evaluate $e_n^TP_nG(-s)z$.
In \S 3, some results obtained with the use of the ${\cal O}(n\log\, n)$ Pascal matrix-vector multiplication
algorithm coupled with a Horner-type scheme
in the computation of B\'ezier curves of degree $n-1$ are then presented.
When $n$ increases, the fast matrix-vector multiplication  becomes unstable, as well as the evaluation of $B(s)$
for $s$ close to 1. Then, the matrix-vector multiplication must be done in a more precise way,
and a stabilizing procedure to  $B(s)$-evaluation should be attempted, e.g. by dividing
the evaluation process in two: first, compute $B(s)$ for $s\in [0,1/2)$; then, compute $B_r(s)$
for $s\in [0,1/2]$, where $B_r(s)$ is the reverse B\'ezier curve, that is,
$B_r(s)$ is the B\'ezier curve defined from the points $P_{n-1}$, ..., $P_1$ and $P_0$,
in this order. For $n>32$, even the 2-steps $B(s)$-evaluation  becomes unstable, yielding
incorrect values for $s$ close to 1/2.
It is when we introduce an affine transform of the vectors of coordinates in order to improve the evaluation.


\section{Pascal matrix-vector multiplication}

In this section some algorithms of Pascal-type matrix-vector
multiplication are discussed from their time complexity.
These algorithms are founded on algebraic properties
of these matrices which are also presented in the following.

\subsection{${\cal O}(n^2)$ arithmetic operations methods}

Algorithms of matrix-vector multiplications with $P_n$,
$P_nG_n(t)$ and $B_n^e(t)$, respectively, are here presented in the form
of MATLAB functions,
all of them demanding ${\cal O}(n^2)$ arithmetic operations.

We begin observing that
$$\left(\begin{array}{rrrr} 1 & 0&0&0\\ 1 & 1 & 0 & 0 \\ 1 & 2 & 1 & 0\\
1 & 3 & 3 & 1 \end{array}\right) = \left(\begin{array}{rrrr} 1 & 0&0&0\\ 0 & 1 & 0 & 0 \\ 0 & 1 & 1 & 0\\
0 & 1 & 2 & 1 \end{array}\right) . \left(\begin{array}{rrrr} 1 & 0&0&0\\ 1 & 1 & 0 & 0 \\ 0 & 1 & 1 & 0\\
0 & 0 & 1 & 1 \end{array}\right),$$
that is
$$P_4 = \left(\begin{array}{cc} 1 & 0 \\ 0 & P_3 \end{array}\right) E_1=
\left(\begin{array}{cc} I_2 & 0\\ 0 & P_2 \end{array}\right)  E_2.E_1= E_3.E_2.E_1.$$
Since for $1<j<i\le n$
$$ e_{i-1}^T P_{n-1} E_{n-1} e_{j-1} = e_{i-1}^T P_{n-1} \left( e_{j-1}+e_j\right)=
e_{i-1}^T P_{n-1} e_{j-1}+ e_{i-1}^T P_{n-1} e_j =$$
$$=\binom{i-2}{j-2} + \binom{i-1}{j-1} = \binom{i-1}{j-1}=\left( P_n\right)_{ij},$$
we conclude that
$$P_n = \left(\begin{array}{cc} 1 & 0 \\ 0 & P_{n-1} \end{array}\right)
\left(\begin{array}{cc} 1 & 0 \\ e^{(n-1)}_1 & E_{n-1} \end{array}\right).$$
Therefore, it has just proved by induction the following statement:

\begin{prop}Let $P_n$ be the $n\times n$ lower triangular Pascal matrix. Then
$P_n = E_{n-1}. \dots . E_1$ where, for all $k\in \{ 1,...,n-1\}$, $E_k = I + e_{k+1}e_k^T + ...
+ e_ne_{n-1}^T$.
\label{paspas}
\end{prop}

The first algorithm is described as a MATLAB function in the following:

\begin{verbatim}
function x = pascal_product(x)
%PASCAL_PRODUCT Multiply a vector x by the lower triangular
%Pascal matrix
n= length(x);
x = x(:);
for k = 2:n
    for s = n:-1:k
        x(s) = x(s) + x(s-1);
    end
end
\end{verbatim}

\

Notice that from the above factorization we also conclude that
$$P_n=\left(\begin{array}{cc} I_k & 0 \\ 0 & P_{n-k} \end{array}\right)
\left(\begin{array}{cc} P_k & 0 \\ Z & W\end{array}\right),$$
where
$$\left(\begin{array}{cc} Z & W \end{array}\right)=
\left(\begin{array}{ccccccc} \binom{k}{0} & \binom{k}{1} &  ... &\binom{k}{k}& 0 & ... & 0\\
0 & \binom{k}{0} & \binom{k}{1} &  ... & \binom{k}{k}& ... & 0\\
\vdots & \ddots & \ddots & \ddots & ... & \ddots & \vdots \\
0 & ... & 0 & \binom{k}{0} & \binom{k}{1} & ... &\binom{k}{k}\end{array}\right).$$

A similar factorization has $P_nG_n(t)$:
$$
\left(\begin{array}{ccccc}
1 & 0 & 0 &   ... & 0 \\
1 & \binom{1}{1}t & 0  & ... & 0 \\
1 & \binom{2}{1}t & \binom{2}{2}t^2  & ... & 0\\
\vdots & \vdots & \vdots & \ddots   & \vdots\\
1 & \binom{n-1}{1}t & \binom{n-1}{2}t^2 &  ... & \binom{n-1}{n-1}t^{n-1}
\end{array}\right)=
\left(\begin{array}{ccccc} 1 & 0 & ... & 0 & 0 \\
                           0 & 1 & ...&0 & 0 \\
                           \vdots  & \vdots & \ddots & \ddots  & \vdots\\
                           0 & 0 & ... &1 & 0 \\
                           0 &  0 & ...    & 1 & t
\end{array}\right)...
\left(\begin{array}{ccccc} 1 & 0 & 0 &...  & 0 \\
                           1 & t & 0& ...& 0 \\
                           0 & 1 & t & ... & 0 \\
                           \vdots  & \vdots & \ddots & \ddots  & \vdots\\
                           0      &  0 & ...    & 1 & t
\end{array}\right),
$$
that is
\begin{prop}Let $P_n$ be the $n\times n$ lower triangular Pascal matrix
and $G(t) = diag(1,t,t^2,...,t^{n-1})$. Then
$P_nG_n(t) = E_{n-1}(t)...E_1(t)$ where, for $1\le k\le n-1$,
$$E_k(t)=e_1e_1^T + ... + e_ke_k^T + e_{k+1}[e_k+te_{k+1}]^T+...+e_n[e_{n-1}+te_n]^T.$$
\label{pasg}
\end{prop}

The second algorithm, which utilizes Proposition~\ref{pasg}, is displayed just below:

\begin{verbatim}
function x = pascal_g_product(x,t)
%PASCAL_G_PRODUCT Multiply a vector x by PG(t) where
%  P is the lower triangular Pascal matrix and
%  G(t) = diag(1,t,t^2,...,t^{n-1})
n= length(x);
x = x(:);
for k = 2:n
    for s = n:-1:k
        x(s) = x(s-1) + t*x(s);
    end
end

\end{verbatim}

By using the notation of \cite{cav}, a Bernstein matrix can be described as
$$ B^e_n(t) = P_n[1-t] G_n(t),$$
where $\left( P_n[t]\right)_{ij} =  \left( P_n\right)_{ij} t^{i-j}$, $1\le j\le i\le n$.
Therefore, it is not difficult to conclude that
a Bernstein matrix has the following bidiagonal factorization:
$$
\left(\begin{array}{ccccc}
1 & 0 & 0 &   ... & 0 \\
1-t & \binom{1}{1}t & 0  & ... & 0 \\
(1-t)^2 & \binom{2}{1}(1-t)t & \binom{2}{2}t^2  & ... & 0\\
\vdots & \vdots & \vdots & \ddots   & 0\\
(1-t)^{n-1} & \binom{n-1}{1}(1-t)^{n-2}t & \binom{n-1}{2}(1-t)^{n-3}t^2 &  ... & \binom{n-1}{n-1}t^{n-1}
\end{array}\right)=
$$
$$=
\left(\begin{array}{ccccc} 1 & 0 & ... & 0 & 0 \\
                           0 & 1 & ...&0 & 0 \\
                           \vdots  & \vdots & \ddots & \ddots  & \\
                           0 & 0 & ... &1 & 0 \\
                           0 &  0 & ...    & 1-t & t
\end{array}\right)...
\left(\begin{array}{ccccc} 1 & 0 & 0 &...  & 0 \\
                           1-t & t & 0& ...& 0 \\
                           0 & 1-t & t & ... & 0 \\
                           \vdots  & \vdots & \ddots & \ddots  & \\
                           0      &  0 & ...    & 1-t & t
\end{array}\right),
$$
that is
\begin{prop} Let $B_n^e(t)$ be a $n\times n$ Bernstein matrix.
Then
$$B_n^e(t)=E_{n-1}^e(t)...E_1^e(t) \mbox{ where, for $1\le k\le n-1$},$$
$$E_k^e(t)=e_1e_1^T + ... + e_ke_k^T + e_{k+1}[(1-t)e_k+te_{k+1}]^T+...+e_n[(1-t)e_{n-1}+te_n]^T.$$
\label{pasb}
\end{prop}
Proposition~\ref{pasb} yields the following algorithm, that is essentially the Casteljau's.

\begin{verbatim}
function x = bernstein_product(x,t)
%BERNSTEIN_PRODUCT Multiply a vector x by a Bernstein matrix
n= length(x);
x = x(:);
t1 = 1-t;
for k = 2:n
    for s = n:-1:k
        x(s) = t1*x(s-1) + t*x(s);
    end
end
\end{verbatim}

Similarly we can conclude that
$$B^e_n(t)=\left(\begin{array}{cc} I_k & 0 \\ 0 & B^e_{n-k}(t) \end{array}\right)
\left(\begin{array}{cc} B^e_k(t) & 0 \\ Z & W\end{array}\right),$$
where
$\left(\begin{array}{cc} Z & W \end{array}\right)$ is the matrix
$$\left(\begin{array}{ccccccc} (1-t)^k & \binom{k}{1}(1-t)^{k-1}t &  ... &\binom{k}{k}t^k& 0 & ... & 0\\
0 & (1-t)^k & \binom{k}{1}(1-t)^{k-1}t &  ... & \binom{k}{k}t^k& ... & 0\\
\vdots & \ddots & \ddots &  ... &  & \ddots & \vdots \\
0 & ... & (1-t)^k  & \binom{k}{1}(1-t)^{k-1}t  &... &   &\binom{k}{k}t^k\end{array}\right).$$

\subsection{Balancing a lower triangular Pascal matrix}

Let $Z_0 = (x_0,y_0)$, $Z_1=(x_1,y_1)$, ..., $Z_{n-1}=(x_{n-1},y_{n-1})$
be $n$ points of $\mathbb{R}^2$.
Let $x = ( x_0 \, x_1 \, \ldots \, x_{n-1} )^T$
and $y = (y_0 \, y_1 \, \ldots \, y_{n-1} )^T$.
For $s \in [0,1]$, let $B^e_n(s)$ be the $n\times n$ Bernstein matrix.
A spectral decomposition of $B^e_n(s)$ is
$B^e_n(s) = P_n G_n(s) P_n^{-1}$, where $P_n$ is the $n\times n$
lower triangular Pascal matrix
and $G_n(s)=diag(1,s,...,s^{n-1})$ (\cite{aceto}). Since $P_n^{-1}=G_n(-1)P_nG_n(-1)$ (\cite{aceto}),
we conclude that the coordinates of the B\'ezier curve $B(s)$, $s\in [0,1]$,
defined from $Z_0$, ..., $Z_{n-1}$ (in this order)
is given by
\begin{equation}
x(s) = e_n^TP_nG_n(-s)P_nG_n(-1)x,
\mbox{ } y(s) = e_n^TP_nG_n(-s)P_nG_n(-1)y. \label{pasc}
\end{equation}
The computation of a matrix-vector multiplication with the Pascal
matrix $P_n$ can be carried out in ${\cal O}(n\log n)$ operations \cite{wang}.
Hence, for each $s\in [0,1]$,
the computation of $x(s)$ and $y(s)$ also requires ${\cal O}(n\log n)$ operations.
However, arithmetic operations with Pascal matrices
are very unstable because of the various magnitudes their entries have.
In this section, a preconditioning technique is proposed to deal with
this instability.
This technique, which is introduced in \cite{wang},
is based on the factorization
\begin{equation}
P_n = D_n(t)\, T_n(t)\, D_n(t)^{-1}, \label{dtd}
\end{equation}
where
$D_n(t)=diag(0!, \dfrac{1!}{t}, \ldots,\dfrac{(n-1)!}{t^{n-1}})$,
and
$$T_n(t) =
\left(\begin{array}{ccccc} 1 & 0 & 0 & \ldots & 0 \\ \dfrac{t}{1!}
& 1 & 0 & \ldots & 0 \\ \frac{t^2}{2!} & \frac{t}{1!} & 1 & \ldots & 0\\
\vdots & \vdots & \vdots & \ddots & \vdots \\
\frac{t^{n-1}}{(n-1)!} & \frac{t^{n-2}}{(n-2)!} &
\frac{t^{n-3}}{(n-3)!} & \ldots & 1 \end{array} \right).$$
Note that for $0<t\le 1$ the entries of the Toeplitz matrix $T_n(t)$
vary from 1 to $\frac{t^{n-1}}{(n-1)!}$,
that is to say $T_n(t)$ is ill-conditioned.
In \cite{wang} it was found that a good value for $t$ is $t_1= \dfrac{n-1}{e}$.
Here we propose a more accurate value for $t$, for almost all values of $n$:
$t_2 = \sqrt[n-2]{\dfrac{(n-1)!}{k+1}}$,
where $k= \lfloor \sqrt[n-1]{(n-1)!} \rfloor$.

In order not to cause great instability in matrix-vector
multiplication with Pascal matrices, we will balance
these matrices by using the decomposition given in (\ref{dtd}).
Since the entries of the matrix $T_n(t)$ are of the form
$$
f(m) = \dfrac{t^m}{m!},
$$
where $m=0,1, \ldots, n-1$,
we would like to find out a value $t$
for which $\max \, f/\min \, f$
is the least possible.
If this optimum $t$ exists, it will bring those entries to be
the closest in magnitude to one another.
The following lemma resumes basic facts about the monotony of $f$.

\begin{lem} Let $t\in (0,\infty)$. Consider the function $f: \N \to [0,\infty)$
defined by $f(m) = \dfrac{t^m}{m!}$. Then
\begin{enumerate}
\item[(i)] $f$ is a nondecreasing function for integers $m$ such that $m
\leqslant \lfloor t \rfloor \leqslant t$;
\item[(ii)] $f$ is a nonincreasing function for integers $m$ such that $t \leqslant \lceil t
\rceil \leqslant m$.
\end{enumerate}
\end{lem}

\begin{pf}
\begin{enumerate}
\item[(i)] $f(m-1) = \dfrac{t^{m-1}}{(m-1)!} =
\dfrac{m}{m}\dfrac{t^{m-1}}{(m-1)!} \leqslant \dfrac{t^m}{m!} =
f(m)$, since $m \leqslant t$.
\item[(ii)] $f(m-1) = \dfrac{t^{m-1}}{(m-1)!} = \dfrac{m}{m}\dfrac{t^{m-1}}{(m-1)!}
\geqslant \dfrac{t^{m}}{m!} = f(m)$, for $m \geqslant t$.
\end{enumerate}\qed
\end{pf}

Therefore, if the optimum value of $t$ exists it would belong to the interval $[1,n-1)$. Now,
$\max \, f/\min \, f=
\dfrac{t^{\lfloor t \rfloor}}{(\lfloor t \rfloor)!}$, if $0\le m\le t$;
if $t<m\le n-1$, then $\max \, f/\min \, f=
\dfrac{t^{\lceil t \rceil} (n-1)!}{t^{n-1}(\lceil t \rceil)!}$.
That is, for the optimum value
$$\min_{t\in [1,n-1)} \max \{ \dfrac{t^{\lfloor t \rfloor}}{(\lfloor t \rfloor)!},
\dfrac{t^{\lceil t \rceil} (n-1)!}{t^{n-1}(\lceil t \rceil)!} \}
$$
is achieved. The following lemma is about the monotony of $f_1$ and $f_2$.
The proof is an exercise of Calculus and will be omitted here.

\begin{lem} Let $n$ be an integer greater than 1.
Let $f_1$ and $f_2$ be two functions defined at
the subset $(0,n-1)$ of real numbers such that
$f_1(t) = \dfrac{t^{\lfloor t \rfloor}}{(\lfloor t \rfloor)!}$
and $f_2(t) = \dfrac{t^{\lceil t \rceil} (n-1)!}{t^{n-1}(\lceil t \rceil)!}$,
respectively. Then
\begin{enumerate}
\item[(i)] $f_1$ is an increasing function;
\item[(ii)] $f_2$ is a nonincreasing function.
\end{enumerate}
Moreover, $f_1$ is continuous and $f_2$ is continuous but at integer values.
\end{lem}

\begin{rem}
Therefore, if exists $\tilde t\in (0,n-1)$ such that $f_1(\tilde t)=f_2(\tilde t)$, this is the optimum value.
\end{rem}

Now, we can calculate the optimum value when it exists.

\begin{lem}
If exists $\tilde{t} \in (k,k+1)$ such that
$f_1(\tilde{t}) = f_2(\tilde{t})$, then
$$\tilde{t} = \sqrt[n-2]{\dfrac{(n-1)!}{k+1}}.$$
\end{lem}
\begin{pf}
Let $\tilde{t}$ be such that $f_1(\tilde{t}) = f_2(\tilde{t})$. Since
$\tilde{t} \in (k,k+1)$, it follows that
$$
\dfrac{\tilde{t}^k}{k!} =
\dfrac{\tilde{t}^{k+1}\, (n-1)!}{(k+1)!\, \tilde{t}^{n-1}} \Rightarrow
\dfrac{k+1}{(n-1)!} = \dfrac{1}{\tilde{t}^{n-2}} \Rightarrow
\tilde{t} = \sqrt[n-2]{\dfrac{(n-1)!}{k+1}}.
$$\qed
\end{pf}

\begin{lem}
If there is a $\tilde{t} \in (k,k+1)$ such that
$f_1(\tilde{t}) = f_2(\tilde{t})$, then
$k = \lfloor \sqrt[n-1]{(n-1)!} \rfloor$.
\end{lem}
\begin{pf}
From the former lemma, if $k < \tilde{t} < k+1$
then ${\tilde{t}}^{n-2} = \dfrac{(n-1)!}{k+1}$.
Hence,
$$
k^{n-2} < \dfrac{(n-1)!}{k+1} < (k+1)^{n-2}.
$$
But this yields
$$
(n-1)! < (k+1)^{n-1} \ \mbox{ and} \ (n-1)!>(k+1)k^{n-2}>k^{n-1}.
$$
Therefore,
$$k = \max\{m \in \Z^+ | \, m^{n-1}<(n-1)!\} = max\{m
\in \Z^+ | \, m<\sqrt[n-1]{(n-1)!}\}
$$
That is,
$k=\lfloor \sqrt[n-1]{(n-1)!} \rfloor.$\qed
\end{pf}

\

\begin{rem} Note that, for $n=3$, $f_1(1)=1$,
$f_2(1)=2$, $f_2(1^+)=1$. So, there is not the optimum value because
$f_2(1) > f_1(1) \ge f_2(1^+)$ and a number greater than 1 but very close to it can be taken as a good value
to balance the $3\times 3$ Pascal matrix.
There are several positive integers $n$, $n>3$, such that for
 some integer $k$, $1\le k\le n-2$,
$f_2(k)> f_1(k) \ge f_2(k^+)$. And
this happens if and only if
$$ \frac{k^k(n-1)!}{k!\, k^{n-1}} > \frac{k^k}{k!} \ge
   \frac{k^{k+1}(n-1)!}{(k+1)!\, k^{n-1}},$$
that is,
$$k^{n-1} < (n-1)! \mbox{ and } k^{n-1} + k^{n-2} \ge (n-1)!.$$
\end{rem}

\

\begin{prop} If there is some integer number $n$, $n>3$,
such that $k^{n-1} < (n-1)!$ and $k^{n-1} + k^{n-2} \ge (n-1)!$
for some integer $k$, $1\le k \le n-2$, then
$$\frac{n-1}{e}<k<\frac{n-1}{2}.
$$
\end{prop}
\begin{pf}
If $k\ge (n-1)/2$, then
$$k^{n-1}\ge (\dfrac{n-1}{2})^{n-1}=\dfrac{(n-1)^{n-1}}{2^{n-1}}.$$
Since, for $n\ge 1$, $n^n\ge (n!)^2$, we have
$$\dfrac{(n-1)^{n-1}}{2^{n-1}} \ge \dfrac{((n-1)!)^2}{2^{n-1}}.$$
Now, for $n>3$, $2^{n-1}<(n-1)!$. Hence
$k^{n-1}>(n-1)!$.

On the other hand, $k$ must be greater than $(n-1)/e$, because
$$k\le \frac{n-1}{e} \Longrightarrow
k^{n-1} + k^{n-2} \le (\dfrac{n-1}{e})^{n-1} [1+\dfrac{e}{n-1}]
$$
$$ < \sqrt{2\pi (n-1)} (\dfrac{n-1}{e})^{n-1}<(n-1)!,$$
by the Stirling formula. \qed
\end{pf}

\begin{figure}[!htb]
\centering
\includegraphics[scale=0.7]{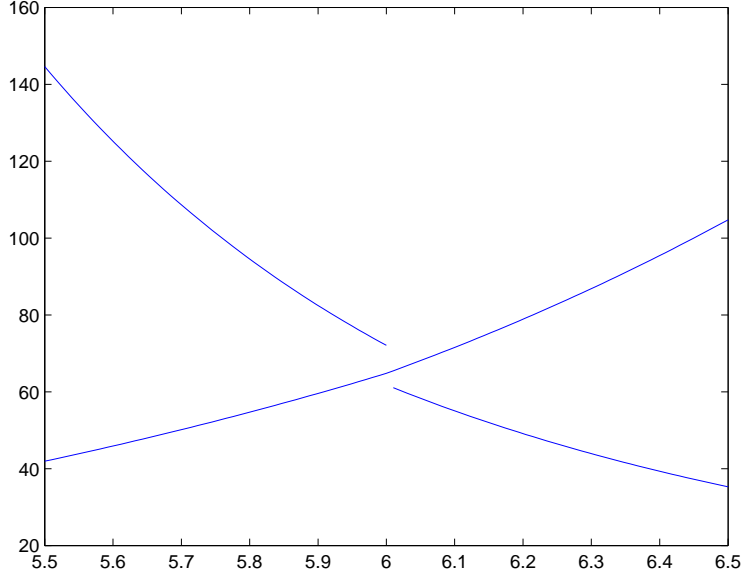}
\caption{$n=15$: $f_1$ and $f_2$ don't intersect}
\end{figure}

\begin{rem} As far as we know, the conjecture about the set of integers
$n$ such that $k^{n-1} < (n-1)! \le k^{n-1} + k^{n-2}$
for some integer $k$, $1\le k \le n-2$, be finite has not been proved yet.
Each integer $n$ less than 10000 belonging to this set
is displayed together with its corresponding integer $k$, $1\le k \le n-2$, in Table~\ref{teste0}.
\end{rem}

\begin{table}[!h]
\centering
\begin{tabular}{|r|r|}
\hline
n&k\\
\hline \hline
15 & 6\\
\hline
39 & 15 \\
\hline
74 & 28 \\
\hline
527& 195\\
\hline
3171 & 1168\\
\hline
5908 & 2175\\
\hline
7036 & 2590\\
\hline
7534 & 3194\\
\hline
7537 & 3401\\
\hline
\end{tabular}

\

\caption{$k^{n-1} \le (n-1)! \mbox{ and } k^{n-1} + k^{n-2} \ge (n-1)!$}
\label{teste0}
\end{table}

For each value of $n$ in the table, we will take
$t = k$ as the optimum value.
For the other values of $n$, $4\le n\le 10000$,
the optimum value is given by
$$
t = \sqrt[n-2]{\dfrac{(n-1)!}{k+1}}
$$
where $k= \lfloor \sqrt[n-1]{(n-1)!} \rfloor$.

In Table~\ref{teste1}
we can compare the accuracy
of the computation of $Pz$ for some values of $n$, where $z$ is the vector
defined by $z_k = (-1)^k$ for $k=1,...,n$:
$w_1=D_n(1).T_n(1).D_n(1)^{-1}z$; 
$w_2=D_n(t).T_n(t).D_n(t)^{-1}z$, with $t=(n-1)/e$;
$w_3=D_n(t).T_n(t).D_n(t)^{-1}z$,
with $t=\sqrt[n-2]{\dfrac{(n-1)!}{k+1}}$
and $k= \lfloor \sqrt[n-1]{(n-1)!} \rfloor$,
or $t=k$ if $(n,k)$ is one of the pairs in Table~\ref{teste0};
$w_4=pascal\_product(z)$.
For the computation of $w_1$, $w_2$ and $w_3$
a fast Toeplitz matrix-vector multiplication was used,
which demanded only ${\cal O}(n\log n)$ operations.
Remark that in exact arithmetic $P_nz=e_1$, the first canonical vector.

\begin{table}[!h]
\begin{center}
\begin{tabular}{|r|r|r|r|c|}
\hline
n&$||P_nz-w_1||_2$ & $||P_nz-w_2||_2$ & $||P_nz-w_3||_2$ & $||P_nz-w_4||_2$\\
\hline \hline
4  & 4.9262e-16  & 4.4790e-16   &  9.6908e-17 &  0\\
\hline
8  & 2.4685e-13  & 9.5720e-15   &  6.2521e-15 &  0\\
\hline
15 & 1.7782e-06 & 5.4617e-13  & 9.3638e-13 &  0\\
\hline
16 & 1.9051e-05 & 9.5411e-13  & 8.4432e-13 &  0\\
\hline
32 & 1.0318e+16 & 2.7147e-07  & 9.5906e-08 &  0\\
\hline
39 & 2.8076e+27 & 1.6654e-05  & 4.4418e-06 &  0\\
\hline
64 & 2.0661e+71 & 1.7187e+04  & 7.2813e+03 &  0\\
\hline
\end{tabular}

\

\caption{Errors in Pascal matrix-vector multiplication}
\label{teste1}
\end{center}
\end{table}

\section{On computing a B\'ezier curve}

In this section we are going to compute a B\'ezier curve.
First, notice that if the control points are translated by a vector $v=(p,q)$,
the B\'ezier curve is also translated by $v$. Another feature of a B\'ezier curve is that
an uniform scaling of the control points yields an uniform scaling of the curve.
Hence, without loss of generality we assume that the coordinates of the control points are all positive
and less than or equal to 1 and, therefore, the $\infty$-norms of $P_nx$ and $P_ny$ are less than or equal
to $2^n$.

\subsection{$B(s)$-evaluation}
Pascal matrix methods compute a B\'ezier curve $B(s)$ of degree $n-1$ via
the decomposition $B^e_n(s) = P_n G_n(-s) P_nG_n(-1)$:
first $z=P_nG_n(-1)x = P_nx_-$  and $w=P_nG_n(-1)y = P_ny_-$ are computed;
then the polynomials $e_n^T P_n G_n(-s)z$ and $e_n^T P_n G_n(-s)w$ are evaluated.
When $n$ is small, e.g. $n=32$, $z$ and $w$ have 2-norms around $10^9$ for $||x||_{\infty}=1$
or $||y||_{\infty}=1$,
and both polynomials could be efficiently evaluated for each $s \in [0,1]$.
The function {\em g\_pascal\_product} does the evaluation, with $n(n-1)/2$ additions and
$n(n-1)/2$ multiplications for each $s$.
A less expensive alternative is to use
a Horner-like scheme that evaluates the polynomial concomitantly with the binomial coefficients.

\

\begin{table}[!b]
\begin{center}
\begin{tabular}{|r|r|r|r|}
\hline
n&$ ||B_P-B_C||_{\infty}$ & $time_P$ & $time_C$\\
\hline \hline
4  & 7.7716e-16  &   0.001810  &  0.001273\\
\hline
8  & 2.8547e-14  &   0.001954  &  0.002164\\
\hline
15 & 9.3585e-11  &   0.002095  &  0.005024\\
\hline
16 & 1.9592e-10  &   0.002238  &  0.005475\\
\hline
24 & 1.2341e-06  &   0.002518  &  0.009927\\
\hline
32 & 0.0190      &   0.003169  &  0.015896\\
\hline
\end{tabular}

\

\caption{Pascal matrix method $\times$ Casteljau's}
\label{first}
\end{center}
\end{table}

The test control points were defined by the MATLAB function $rand(n,2)$,
which returns an $n\times 2$ matrix containing pseudo-random values
drawn from a uniform distribution on the unit interval. We have used the Pascal matrix-vector
multiplication done from the similar Toeplitz matrix $T(t)$, with $t$ found by our procedure
(according to our calculations in the last section), plus the $B(s)$ evaluation given by the Horner-like
scheme cited above. The results were compared with the ones obtained by Casteljau's method.
In table~\ref{first} the $\infty$-norms of the differences of the results obtained by
the Pascal matrix method ($B_P$) and Casteljau's ($B_C$) are listed close by the average time of computation
of all B(s) evaluation ($\Delta\, s = 1/128$). We have taken the smallest time among 10 elapsed times obtained from consecutive executions of the procedure as the average time of computation, all of them computed by the MATLAB's built-in tic/toc functions.

For $n=32$, the $B(s)$-evaluation becomes unstable when $s$ approaches to 1.
To locally and globally improve the evaluation, we have made
a simple procedure which has yielded more precise results, that is
to divide the process of evaluation in two independent steps:
\begin{enumerate}
\item[(a)] evaluate $e_n^T P_n G_n(-s)z$ and $e_n^T P_n G_n(-s)w$ for $0\le s\le 1/2$;
\item[(b)] evaluate $e_n^T P_n G_n(-s)z_r$ and $e_n^T P_n G_n(-s)w_r$ for $1/2> s\ge 0$, which
is equivalent to evaluate $e_n^T P_n G_n(-s)z$ and $e_n^T P_n G_n(-s)w$ for $1/2 < s\le 1$.
\end{enumerate}
The results indicate that the procedure has worked well as far as $n=41$ and
some of those can be seen in table\ref{sec}.

\begin{table}[!hb]
\begin{center}
\begin{tabular}{|r|r|r|r|}
\hline
n&$ ||B_P-B_C||_{\infty}$ & $time_P$ & $time_C$\\
\hline \hline
32  & 2.3113e-07 &   0.005602  &  0.015901\\
\hline
36  & 2.6961e-05 &   0.006106  &  0.019464\\
\hline
39  & 1.3152e-04 &   0.005873  &  0.022314\\
\hline
41  & 4.8668e-04 &   0.006229  &  0.024372\\
\hline
42  &  0.0022 &   0.006387  &  0.025476\\
\hline
48  &  0.1112 &   0.006706  &  0.032366\\
\hline
\end{tabular}

\

\caption{Pascal matrix method with reverse evaluation $\times$ Casteljau's}
\label{sec}
\end{center}
\end{table}

Two remarkable facts arise in table~\ref{sec}: first, the time of computation for $n=39$
is smaller than that for $n=41$;
second, the sudden loss of precision from $n=41$ to $n=42$.
One of the explanations for the first fact is because there is no calculation to find
the optimum value to balance $P_{39}$: 15 is taken to be the optimum value.
The second fact surely has to do with the limitations of the flowing point arithmetic
of our machine, a 32-bits AMD Athlon XP 1700+ (1467 MHz).

\subsection{On conditioning the vectors of coordinates}

One way to overcome this lack of stability is to transform
the vectors $x$ and $y$ of coordinates into a vector very near to $e^T=( 1\, 1\, ...\, 1)^T$.
Since
$$v_k = \frac{1}{m+1} v_{k-1} + \frac{m}{m+1}e$$
is a stationary scheme that
converge to the solution $e$ of $Ix=e$ for any $v_0$,
the idea is to compute the B\'ezier curve $T_m(B)$ from control
points $W_0=T_m(Z_0)$, ..., $W_{n-1}=T_m(Z_{n-1})$, where $T_m(v) = (v+m.e)/(m+1)$, and then
to obtain $B$ by inverse transforming the points of $T(B)$. Note
that
$$T_m(B)(s)=e_n^T P_nG(-s) P_n G(-1) \left(\begin{array}{c} T_m(x) // T_m(y)\end{array}\right).$$
Hence,
$$e_n^T P_nG(-s) P_n G(-1) T_m(x) = e_n^T P_nG(-s) P_n G(-1)(x+m.e)/(m+1)=$$
$$=\frac{1}{m+1}e_n^T P_nG(-s) P_n G(-1)x + \frac{m}{m+1}e_n^T P_nG(-s) P_n G(-1)e =$$
$$=\frac{1}{m+1}x(s) +  \frac{m}{m+1}.$$
Thus, since an analogous result is obtained with $y$,
$x(s) = (m+1) T_m(B)(s) -m$.
Observe that
$$T_m T_{m'} (v) = \left( v + [(m+1)(m'+1)-1].e \right) / [(m+1)(m'+1)] = T_{(m+1)(m'+1)-1}.$$
Our strategy to have better $B(s)$-evaluation is the following:
\begin{enumerate}
\item[(a)] evaluate $e_n^T P_n G_n(-s)z$ and $e_n^T P_n G_n(-s)w$ for $0\le s\le 1/3$;
\item[(b)] evaluate $e_n^T P_n G_n(-s)P_nG_n(-1)T_m(x)$ and $e_n^T P_n G_n(-s)P_nG_n(-1)T_m(y)$ for $1/3< s < 2/3$,
and inverse transform;
\item[(c)] evaluate $e_n^T P_n G_n(-s)z_r$ and $e_n^T P_n G_n(-s)w_r$ for $1/3\ge s\ge 0$.
\end{enumerate}

For the experiments, the coordinates of the $n$ control points were defined
from the command $A=rand(n,2)$, followed by the normalization $A=A/norm(A)$.
$T_m\circ T_m$, where $m+1=32768$ ($2^{15}$), were applied to $A$ when $n=42,48,54$;
$T_{m'}\circ T_m\circ T_m$, where $m'+1 = 1024$, when $n=59$;
$T_{m''}\circ T_{m'}\circ T_m\circ T_m$, where $m''+1=4$, when $n=64$.
The results are in table\ref{third}.

\begin{table}[!h]
\begin{center}
\begin{tabular}{|r|r|r|r|}
\hline
n&$ ||B_P-B_C||_{\infty}$ & $time_P$ & $time_C$\\
\hline \hline
42  & 8.3290e-07 &   0.013525  &  0.025363\\
\hline
48  & 1.7620e-06 &   0.015566  &  0.032190\\
\hline
54  & 2.3903e-04 &   0.017786  &  0.039937\\
\hline
59  & 9.9235e-04 &  0.019782  &  0.047320\\
\hline
64  &  0.0048 &  0.022044  &  0.055178\\
\hline
\end{tabular}

\

\caption{Pascal matrix method with piecewise evaluation $\times$ Casteljau's}
\label{third}
\end{center}
\end{table}

\section{Conclusions}

We have presented results obtained from some methods to compute a B\'ezier curve of degree $n-1$,
for various values of $n$. They
were created from a description of the curve that involves matrix-vector multiplications with
the $n\times n$ lower triangular Pascal matrix $P_n$, which are here called Pascal matrix methods. With this in mind we have introduced two algorithms: one, which only demands $n(n-1)/2$ additions, is very precise and it is based on the fact that $P_n$ is a product of bidiagonal matrices with 0 and 1; the other, which demands ${\cal O}(n\, \log\, n)$ algebraic operations, depends on a positive
real value in order to minimize the magnitudes of the entries of $P_n$ when considered as a scaled Toeplitz matrix.
We have seen that there is a function that relates $n$ to that optimum value except for the integers belonging to a certain set, which we have not yet known if it is finite or not. Once the matrix-vector multiplication done, a polynomial evaluation
should have been carried out for various values $s\in [0,1]$, which has become unstable as $s$ approaches 1/2.
From the set of experiments presented here we have seen that the combination of Pascal matrix-vector multiplication plus polynomial evaluation has converged to the results obtained from Casteljau's by adopting some strategies, which vary according to the magnitude of $n$. And even so, they
are more effective concerning time of computation than Casteljau's, at least for $n\le60$.

\end{document}